\documentclass[12pt]{article}
\usepackage{amsmath,amsfonts,amsthm,amssymb,color}
\usepackage{fancybox}
\usepackage{framed}
\usepackage{tikz, float, xcolor, color, graphicx}
\usetikzlibrary{calc,intersections}
\usetikzlibrary{arrows}

\newtheorem{theorem}{Theorem}[section]

\theoremstyle{definition}

\newenvironment{fminipage}%
  {\begin{Sbox}\begin{minipage}}%
  {\end{minipage}\end{Sbox}\fbox{\TheSbox}}

\newcommand{\E}{\mbox{{\bf E}}}

\title{Dissecting the Equilateral Triangle into Non-Congruent Equilateral Triangles}
\author{
Timothy Chu
\\MIT \\timchu@mit.edu
}

\begin{document}
\maketitle
\begin{abstract}
In this paper, we show that an equilateral triangle cannot be dissected into finitely many smaller equilateral triangles, no two of which share two vertices. We do this without the use of Electrical Networks.
\end{abstract}
\section{Introduction}
The problem of dissecting a square into finitely many smaller squares, no two of which are the same size, is a famous problem in mathematics dating back to Dehn in 1903. This problem, also known amusingly as 'squaring the square', was approached by famous combinatorialists in the 30's, most notably Brooks, Smith, Stone, and Tutte, who invented a technique involving electrical networks that was used to construct non-trivial rectangles which could be tiled by squares in \cite{BrooksSST40}. For many years, it was thought that it was impossible to 'Square the square' until R. Sprague constructed the first example in 1939, in \cite{S39}.

In this paper, we examine the related problem: can you dissect an equilateral triangle into finitely many equilateral triangles, no two of which are the same size? This related dissection problem remained open until Tutte demonstrated the impossibility of such a dissection in 1947, by utililzing the machinery of topological dualities and electrical networks with modified Kirchoff laws in \cite{Tutte48}.

In this paper, we present two simple and novel proofs of Tutte's result on equilateral triangles, both of which do not use any machinery from electrical networks. The advantage of these approaches are that they are more elementary,  and can perhaps shed a different type of insight into dissection problems.

\section{Two Proofs}
We will prove a slightly stronger statement, which Tutte also proved in his 1947 paper.

\begin{theorem} An equilateral triangle cannot be tiled with smaller equilateral triangles, no two of which have two vertices in common \end{theorem}

\subsection{First Proof: An Extremal Argument}

Suppose there exists a dissection of the equilateral triangle such that no two share two vertices. Consider the smallest instance of the following, where smallest is defined by the smallest side of the larger triangle:

\begin{figure}[H]
\scalebox{1}{
\begin{tikzpicture}

[point/.style = {draw, circle,  fill = black, inner sep = 1pt},]

\def \x{ 0.9 }
\def \A {(-2, 0)}
\def \B {(0, 2*1.732)} 
\def \C {(2, 0)}
\def \D {(2 - \x, \x*1.732)}
\def \E {(2 + \x, \x * 1.732)}
\draw \A -- \B -- \C -- \A;
\draw \C -- \E -- \D;
\draw (2, 0) -- (3,0);

\node[below] at \A {$A$};
\node[above] at \B {$B$};
\node[below] at \C {$C$};
\node [above] at \D {$D$};
\node [above] at \E {$E$};

\end{tikzpicture}
}
\caption{The Figure we are Looking For}
\end{figure}
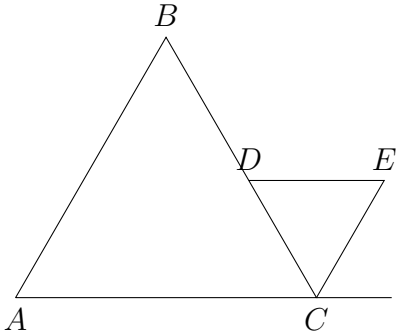

Note this type of figure exists at every corner of the dissected equilateral triangle, and thus the "smallest" instance of this figure must exist.

Now we can use casework, based on the edges in our dissection that have an endpoint at $E$.

\begin{figure}[H]
\scalebox{1}{
\begin{tikzpicture}

\def \x{ 0.9 }
\def \y{1.2} 
\def \A {(-2, 0)}
\def \B {(0, 2*1.732)} 
\def \C {(2, 0)}
\def \D {(2 - \x, \x*1.732)}
\def \E {(2 + \x, \x * 1.732)}
\draw \A -- \B -- \C -- \A;
\draw \C -- \E -- \D;
\draw (2, 0) -- (3,0);

\draw(2, 0) -- (2+\y, \y*1.732); 

\node[below] at \A {$A$};
\node[above] at \B {$B$};
\node[below] at \C {$C$};
\node [above] at \D {$D$};
\node [above] at \E {$E$};

\end{tikzpicture}
}
\caption{Case 1: The line $CE$ extends past $E$}
\end{figure}
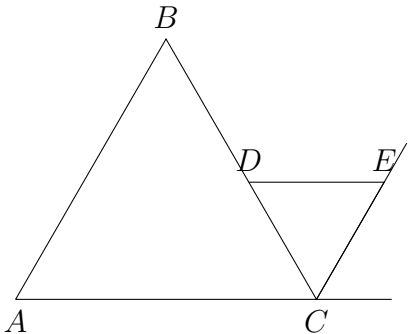

In this case, there must be some triangle with a vertex on $E$ and a side on segment $DE$.

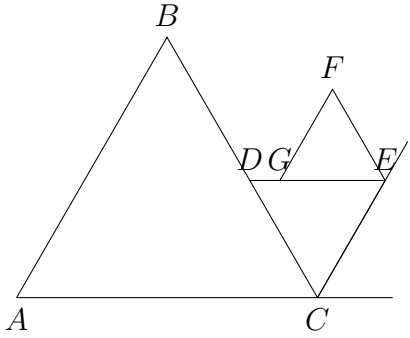
\begin{figure}[H]
\scalebox{1}{
\begin{tikzpicture}

\def \x{ 0.9 }
\def \y{1.2} 
\def \z {0.7}
\def \A {(-2, 0)}
\def \B {(0, 2*1.732)} 
\def \C {(2, 0)}
\def \D {(2 - \x, \x*1.732)}
\def \E {(2 + \x, \x * 1.732)}
\def \F {( 2 + \x - \z, \x * 1.732 + \z * 1.732)}
\def \G {(2+\x - 2 * \z, \x * 1.732)}

\draw \A -- \B -- \C -- \A;
\draw \C -- \E -- \D;
\draw \E -- \F -- \G;
\draw (2, 0) -- (3,0);

\draw(2, 0) -- (2+\y, \y*1.732); 

\node[below] at \A {$A$};
\node[above] at \B {$B$};
\node[below] at \C {$C$};
\node [above] at \D {$D$};
\node [above] at \E {$E$};
\node[above] at \F {$F$};
\node[above] at \G {$G$};

\end{tikzpicture}
}

\caption{Case 1: The line $CE$ extends past $E$}
\end{figure}

This is a contradiction, because figure $CDEFG$ is a smaller version of figure $ABCED$, and we assumed $ABCED$ to be the smallest figure of its kind.

The other two cases, when $DE$ extends beyond $E$, and when a line 60 degrees away from $DE$ and $CE$ intersects $E$, are similar.

\begin{figure}[H]
\scalebox{1}{
\begin{tikzpicture}

\def \x{ 0.9 }
\def \y{1.6} 
\def \z {0.4}
\def \A {(-2, 0)}
\def \B {(0, 2*1.732)} 
\def \C {(2, 0)}
\def \D {(2 - \x, \x*1.732)}
\def \E {(2 + \x, \x * 1.732)}
\def \F {( 2 + \x + \z, \x * 1.732 - \z * 1.732)}
\def \G {(2+\x -  \z, \x  * 1.732 - \z * 1.732)}

\draw \A -- \B -- \C -- \A;
\draw \C -- \E;
\draw \E -- \F -- \G;
\draw (2, 0) -- (3,0);
\draw (2 + \y,  \x*1.732) --  \D; 

\node[below] at \A {$A$};
\node[above] at \B {$B$};
\node[below] at \C {$C$};
\node [above] at \D {$D$};
\node [above] at \E {$E$};
\node[above] at \F {$F$};
\node[above] at \G {$G$};

\end{tikzpicture}
}

\caption{Case 2: The line $DE$ extends beyond $E$.}
\end{figure}
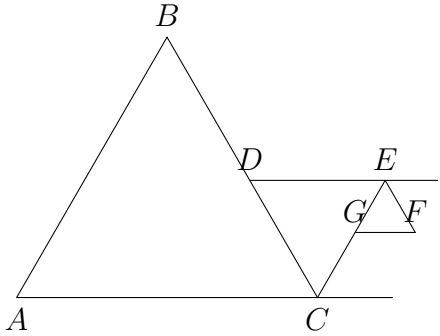

 In this case, there must exist a triangle $FEG$ such that $G$ is on segment $EC$. This implies that $DCEFG$ is a smaller copy of $ABCED$, which is a contradiction.

\begin{figure}[H]
\scalebox{1}{
\begin{tikzpicture}

\def \x{ 0.9 }
\def \y{1.6} 
\def \z {0.5}
\def \A {(-2, 0)}
\def \B {(0, 2*1.732)} 
\def \C {(2, 0)}
\def \D {(2 - \x, \x*1.732)}
\def \E {(2 + \x, \x * 1.732)}
\def \F {(2 - \x + \z, \x * 1.732 + \z * 1.732)}
\def \G {(2- \x + 2 * \z, \x * 1.732)}
\def \X {( 2 + \x + \z, \x * 1.732 - \z * 1.732)}
\def \Y {(2+\x -  \z, \x  * 1.732 + \z * 1.732)}

\draw \A -- \B -- \C -- \A;
\draw \C -- \E -- \D;
\draw \D -- \F -- \G;
\draw \X -- \Y;
\draw (2, 0) -- (3,0); 

\node[below] at \A {$A$};
\node[above] at \B {$B$};
\node[below] at \C {$C$};
\node [above] at \D {$D$};
\node [above] at \E {$E$};
\node[above] at \F {$F$};
\node[above] at \G {$G$};

\end{tikzpicture}
}

\caption{Case 3: There exists a line through $E$ that forms a $60$ degree angle with both $DE$ and $EC$.}
\end{figure}
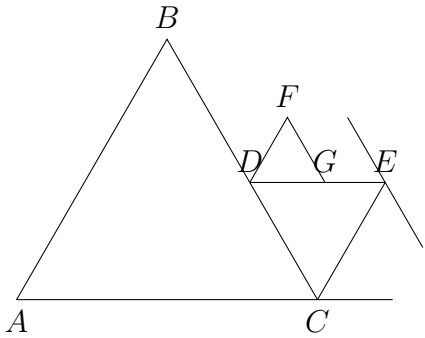

In this last case, there must exist a triangle $FDG$ such that $G$ is on segment $DE$. This implies that $CEDFG$ is a smaller copy of $ABCED$, a contradiction.

This covers all cases, and therefore an equilateral triangle may not be dissected into smaller equilateral triangles, no two of which share two vertices.

\subsection{Second Proof: The Walking Hamster}
The proof in this section is due to Mitchell M. Lee.
\vspace{2 mm}

Suppose that there exists a dissection of the equilateral triangle with smaller equilateral triangles. Now imagine a hamster starting in a corner of the large equilateral triangle, who always takes a 120 degree turn every time he meets a vertex of any smaller equilateral triangle. If no two equilateral triangles have two vertices in common, then a basic induction will prove that the Hamster's path consists of smaller and smaller line segments. It follows that the Hamster's path has no cycles. However, the hamster can also never stop walking; therefore any dissection of the equilateral triangle into smaller equilateral triangles, no two of which share two vertices, must be an infinite dissection.


\bibliographystyle{alpha}

\end{document}